\documentclass[11pt]{article}
\usepackage[margin=1in]{geometry}
\usepackage{amsthm,amsmath,amsfonts,amssymb,mathtools,mathrsfs}
\usepackage[numbers]{natbib}
\numberwithin{equation}{section}
\theoremstyle{plain}
\newtheorem{theorem}{Theorem}[section]
\newtheorem{proposition}[theorem]{Proposition}
\newtheorem{lemma}[theorem]{Lemma}
\newtheorem{corollary}[theorem]{Corollary}
\theoremstyle{remark}

\newcommand{\E}{\mathbb E}
\newcommand{\Pp}{\mathbb P}
\newcommand{\F}{\mathcal F}
\newcommand{\R}{\mathbb R}
\newcommand{\one}{\mathbf 1}

\newcommand{\sinc}{\operatorname{sinc}}
\newcommand{\HS}{\mathrm{HS}}

\usepackage[hidelinks]{hyperref}
\title{Exponential Asymptotics for Brownian Stochastic Integrals with Delayed Information}
\author{Yan Dolinsky\\[3pt]\small Hebrew University, Jerusalem, Israel\\
\small\texttt{yan.dolinsky@Mail.huji.ac.il}}
\date{}
\begin{document}
\maketitle
\begin{abstract}
We study the effect of a small observation delay on the approximation of Brownian stochastic integrals. 
In a regime where the risk-aversion parameter is inversely proportional to the delay, we derive a 
large-deviation-type scaling limit for the Brownian stochastic integral. The limit is expressed 
in terms of an explicit spectral function that captures the effect of the short observation window. 
For nonnegative globally Lipschitz coefficients, the limit is finite when the supremum norm is below $1/2$
and infinite when it is above $1/2$.
As an application, we quantify the exponential hedging error associated with a prescribed delayed-delta strategy.
\end{abstract}


\section{Introduction}
Let $W:=(W_t)_{t\ge0}$ be standard Brownian motion on
$(\Omega,\F,\Pp)$, with its augmented natural filtration
$(\F_t)_{t\ge0}$.
The martingale
representation theorem (see \cite[Section~4]{KaratzasShreve1991}) states
that, for a measurable function $f:\R\to\R$ such that $f(W_1)\in L^2$,
\begin{equation*}
 f(W_1)=u(0,0)+\int_0^1u_x(t,W_t)\,dW_t,
\end{equation*}
where $u=u(t,x)$ solves the backward heat equation
\[
 u_t+\frac12u_{xx}=0,\qquad t<1,
\]
and is given, for $0\le t<1$, by
\[
 u(t,x):=\E\left[f(x+W_{1-t})\right]
 =\frac1{\sqrt{2\pi}}\int_{\R}
 f\bigl(x+\sqrt{1-t}\,y\bigr)e^{-y^2/2}\,dy.
\]
For measurable $f$, the terminal condition is understood in the sense
that $u(t,W_t)\to f(W_1)$ in $L^2$ as $t\uparrow1$.
From the point of view of mathematical finance, the martingale
representation theorem can be interpreted as a perfect replication result,
and $(u_x(t,W_t))_{t\in[0,1)}$ is the corresponding delta-hedging strategy.

We consider a market in which investment decisions are based on delayed
price observations, while trades are executed at the current price. Such
a delay may arise from the time needed to receive and process market information.

For $h>0$, set $\tau_h(t):=(t-h)^+=\max\{t-h,0\}$. The corresponding delayed-delta strategy is $u_x(\tau_h(t),W_{\tau_h(t)})$.
We study the limit as $h\downarrow 0$ of the expression 
\begin{equation}\label{1.1}
h\log\E\left[e^{\frac{1}{h}\left(f(W_1)-\int_{0}^1 u_x\left(\tau_h(t), W_{\tau_h(t)}\right)dW_t\right)}\right].
\end{equation}
This quantity measures the exponential hedging error of the delayed-delta strategy
for a European contingent claim with payoff $f(W_1)$, under (high)
risk aversion $1/h$.

The utility-indifference pricing problem with delayed information was studied in
\cite{BankDolinsky2021} in the Bachelier model, in the same regime in which
risk aversion grows inversely with the observation delay. The analysis there allows
optimization over admissible delayed-information portfolios and identifies the
limiting indifference price through a volatility-control problem with a
quadratic penalty. In this sense, optimization can yield a better value than
that associated with a prescribed trading strategy, when both are evaluated
under the same exponential-utility criterion. However, the limiting value
characterization in \cite{BankDolinsky2021} does not provide an explicit
family of asymptotically optimal delayed-information hedging portfolios for
general payoffs.

For smooth $f$ with suitable integrability, It\^o's formula gives
\begin{equation}\label{eq:ito-identity}
f(W_1)-\int_0^1 u_x\left(\tau_h(t), W_{\tau_h(t)}\right)dW_t=
u(0,0)+\int_0^1\int_{\tau_h(t)}^t u_{xx}(s,W_s)\,dW_s\,dW_t.
\end{equation}
Hence, studying the term in \eqref{1.1}
leads naturally to studying scaling limits for double stochastic integrals.

A large-deviation-type scaling limit for the discretization error of a
Brownian stochastic integral was established in \cite{Dolinsky2027}. An integral of a pointwise spatial supremum
also appears in that limit.
In contrast to the disjoint time cells in that setting, the observation
windows here overlap, leading to a different spectral function and a
critical upper threshold of $1/2$ rather than $1$.

This problem connects the large-deviation theory of Wiener functionals
with the spectral analysis of Gaussian quadratic forms. While
classical exponential asymptotics for small-noise Brownian paths were
established in \cite{Schilder1966}, the small parameter here enters through the
observation delay and the factor $1/h$ in the exponential moment.
For deterministic coefficients, the delayed double It\^o integral belongs
to the second Wiener chaos, and its exponential moment can be analyzed
through the eigenvalues of an associated Hilbert--Schmidt operator, using
the Gaussian quadratic-form formulas in Shepp's seminal paper \cite{Shepp1966}. The kernel of
this operator is concentrated in a band of width $h$ around the diagonal.
After rescaling, its spectral behavior produces the $\sinc$ term in the
limiting function $\Psi$, echoing Szeg\H{o}-type asymptotics for Toeplitz
forms \cite{GrenanderSzego1958}. Related spectral methods were used in
\cite{BrycDembo1997} to study large deviations of quadratic functionals
of stationary Gaussian processes; here, however, the kernel depends on
the vanishing delay and can have nonconstant coefficients. The case of
state-dependent coefficients requires an additional probabilistic argument:
changes of measure and approximation along deterministic paths give the lower bound,
while Gaussian comparison and block estimates control the upper bound.
The change-of-measure argument is closely related to the variational
approach of \cite{BoueDupuis1998}.

The paper is organized as follows. In Section~\ref{sec:preliminary}, we present a preliminary result that motivates the main results stated in Section~\ref{sec:main}.
Section~\ref{sec:auxiliary} collects the auxiliary lemmas used in the proofs.
Section~\ref{sec:lower} establishes the lower bound by a change of measure and deterministic paths, while Section~\ref{sec:upper} proves the upper bound by block estimates.

\section{Preliminary result}\label{sec:preliminary}
We first consider the double It\^o integral
\[
N_h(b):=\int_0^1\int_{(t-h)^+}^{t}b(s)\,dW_s\,dW_t,
\qquad h>0,
\]
where $b:[0,1]\to[0,\infty)$ is continuous (deterministic).

Let
 $\sinc z:=\sin z/z$ for $z\ne0$ and $\sinc(0):=1$.
For $c\in[0,1/2)$, define
\begin{equation*}
 \Psi(c):=\frac1{4\pi}\int_{\mathbb R}
 \left[-\log(1-2c\sinc z)-2c\sinc z\right]dz.
\end{equation*}

The function $\Psi$ is nonnegative and continuous on $[0,1/2)$,
by the inequality $-\log(1-y)-y\ge0$ for $y<1$ and dominated convergence
(using a constant multiple of $\sinc^2 z$ on compact subintervals).
Throughout, $\|\cdot\|_\infty$ denotes the supremum norm on the
relevant domain, allowing the value $+\infty$.

The following proposition gives the asymptotic behavior of $N_h(b)$.
\begin{proposition}\label{prop:deterministic}
Let $b:[0,1]\to[0,\infty)$ be continuous.
\begin{enumerate}
\item If $2\|b\|_\infty<1$, then $\E\left[\exp\{N_h(b)/h\}\right]<\infty$ for every $h>0$, and
\begin{equation*}
\lim_{h\downarrow0}h\log\E\left[\exp\left(\frac{N_h(b)}h\right)\right]
=\int_0^1\Psi(b(t))\,dt.
\end{equation*}
\item If $2\|b\|_\infty>1$, then
\[
\E\left[\exp\left(\frac{N_h(b)}h\right)\right]=+\infty
\qquad\text{for every sufficiently small }h>0.
\]
\end{enumerate}
\end{proposition}

The proof uses the Gaussian quadratic-form identities of \cite{Shepp1966}.

\begin{proof}
Let $\mathcal H:=L^2(0,1)$ and let $A_{h,b}$ be the integral operator on $\mathcal H$
with symmetric kernel
\begin{equation*}
K_{h,b}(s,t):=\frac1h b(\min\{s,t\})\one_{\{|s-t|<h\}},
\qquad s,t\in[0,1].
\end{equation*}
Since $b$ is bounded, $A_{h,b}$ is a compact self-adjoint
Hilbert--Schmidt operator. Write $(\lambda_{h,j})_{j\ge1}$ for its
nonzero eigenvalues, repeated with multiplicity, so that
$\sum_j\lambda_{h,j}^2<\infty$. Let $\lambda_{\max}(A_{h,b})$
denote the top of the spectrum.
From formulas (18.1)--(18.2) in \cite{Shepp1966} we obtain
\begin{equation}\label{eq:sheppformula}
\log\E\left[e^{N_h(b)/h}\right]
=\frac12\operatorname{Tr}\phi(A_{h,b})\in[0,+\infty],
\end{equation}
where
\[
\phi(x):=\begin{cases}
-\log(1-x)-x,&x<1,\\
+\infty,&x\ge1.
\end{cases}
\]
Here the trace is interpreted in the extended sense, as the sum of the
nonnegative spectral contributions.
In particular, the exponential moment is finite if and only if
$\lambda_{\max}(A_{h,b})<1$.

Put $M:=\|b\|_\infty$. Since
$|K_{h,b}(s,t)|\le M h^{-1}\one_{\{|s-t|<h\}}$,
Schur's test \cite[Theorem~5.2]{HalmosSunder1978} yields
\begin{equation}\label{eq:det-schur}
\|A_{h,b}\|_{\mathrm{op}}\le 2M.
\end{equation}
Thus, if $2M<1$, the exponential moment is finite for every $h>0$.

Next, since $A_{h,b}$ is Hilbert--Schmidt, $A_{h,b}^m$ is
trace class for every integer $m\ge2$. Expanding the kernel of
its $m$-fold product and taking the trace gives the cyclic formula
\[
\operatorname{Tr}(A_{h,b}^m)
=h^{-m}\int_{[0,1]^m}\prod_{j=1}^m
b(t_j\wedge t_{j+1})\,
\one_{\{|t_j-t_{j+1}|<h\}}\,dt_1\cdots dt_m,
\qquad t_{m+1}=t_1.
\]
Put $k(y):=\one_{(-1,1)}(y)$ and substitute
$t_1:=t$, $t_j:=t+hy_j$ for $j=2,\ldots,m$, with
$y_1=y_{m+1}=0$. For this change of variables, extend $b$
continuously and boundedly outside $[0,1]$; the domain indicators
below ensure that this does not change the integral.
The Jacobian $h^{m-1}$ cancels all but
one factor of $h^{-m}$, so
\[
\begin{aligned}
h\operatorname{Tr}(A_{h,b}^m)
={}&\int_0^1\int_{\mathbb R^{m-1}}
\prod_{j=1}^m k(y_j-y_{j+1})
\prod_{j=1}^m b\bigl(t+h\min\{y_j,y_{j+1}\}\bigr)\\[-2pt]
&\hspace{30mm}\times
\prod_{j=2}^m\one_{\{0\le t+hy_j\le1\}}
\,dy_2\cdots dy_m\,dt.
\end{aligned}
\]
Since $y_1=y_{m+1}=0$ and the factors $k(y_j-y_{j+1})$
force $|y_j-y_{j+1}|<1$, all $y_j$ are bounded by $m$ on the
support of the integrand. For $t\in(0,1)$, the arguments of $b$
therefore converge to $t$ and the domain indicators converge to $1$.
Boundedness and continuity of $b$ now allow dominated convergence:
\[
\lim_{h\downarrow0}h\operatorname{Tr}(A_{h,b}^m)
=\int_0^1b(t)^m\,dt
\int_{\mathbb R^{m-1}}\prod_{j=1}^m k(y_j-y_{j+1})
\,dy_2\cdots dy_m
=(k^{*m})(0)\int_0^1 b(t)^m\,dt.
\]
The remaining integral is $(k^{*m})(0)$, the value at zero of
the $m$-fold convolution of the even function $k$.
Since $\widehat k(z)=\int_{-1}^1 e^{-izy}\,dy=2\sinc z$,
the convolution theorem and Fourier inversion (see \cite{SteinShakarchi2003}) give
$(k^{*m})(0)=(2\pi)^{-1}\int_{\mathbb R}(2\sinc z)^m\,dz$;
the integral converges for $m\ge2$. Hence
\begin{equation}\label{eq:det-powers}
\lim_{h\downarrow0}h\operatorname{Tr}(A_{h,b}^m)
=\frac1{2\pi}\int_0^1\!\int_{\mathbb R}
(2b(t)\sinc z)^m\,dz\,dt,\qquad m\ge2.
\end{equation}

Suppose $2M<1$. By \eqref{eq:det-schur}, the power series
$\phi(x)=\sum_{m=2}^\infty x^m/m$ converges uniformly on the
spectrum of $A_{h,b}$, and hence
\[
h\log\E\left[e^{N_h(b)/h}\right]
=\frac h2\sum_{m=2}^\infty\frac{\operatorname{Tr}(A_{h,b}^m)}m.
\]
Moreover, for every $m\ge2$,
\[
\left|h\operatorname{Tr}(A_{h,b}^m)\right|
\le (2M)^{m-2}h\operatorname{Tr}(A_{h,b}^2)
\le 2M^2(2M)^{m-2}.
\]
We can therefore pass to the limit term by term using
\eqref{eq:det-powers}. Summing the series under the integral yields
\[
\lim_{h\downarrow0}h\log\E\left[e^{N_h(b)/h}\right]
=\frac1{4\pi}\int_0^1\!\int_{\mathbb R}
\phi(2b(t)\sinc z)\,dz\,dt
=\int_0^1\Psi(b(t))\,dt,
\]
as claimed.

If $2M>1$, continuity of $b$ allows us to choose
$g\in C_c^\infty(0,1)$ with $\|g\|_2=1$ and
$2\int_0^1b(t)g(t)^2\,dt>1$. By the approximate-identity limit,
\[
\lim_{h\downarrow0}\langle g,A_{h,b}g\rangle
=2\int_0^1b(t)g(t)^2\,dt>1.
\]
Hence $\lambda_{\max}(A_{h,b})>1$ for all sufficiently small $h$,
and the finiteness criterion following \eqref{eq:sheppformula}
implies that the exponential moment is infinite for such $h$.
\end{proof}

We conclude this section by noting that the preceding arguments do not cover the critical case $2\|b\|_\infty=1$, which we leave for future research.
\section{Main results}\label{sec:main}
The following theorem is the main result of the paper. 

\begin{theorem}\label{thm:nonlinear}
Let $a:[0,1]\times\R\to[0,\infty)$ be globally Lipschitz, and put
\begin{equation*}
N_h(a):=\int_0^1\int_{\tau_h(t)}^t a(s,W_s)\,dW_s\,dW_t.
\end{equation*}
\begin{enumerate}
\item If $\|a\|_\infty<1/2$, then
\begin{equation}\label{eq:nonlinear-limit}
\lim_{h\downarrow0}h\log\E\left[\exp\left(N_h(a)/h\right)\right]
= \int_0^1\sup_{x\in\R}\Psi(a(t,x))\,dt.
\end{equation}
\item If $\|a\|_\infty>1/2$, then
\begin{equation}\label{eq:nonlinear-supercritical}
\lim_{h\downarrow0}h\log\E\left[\exp\left(N_h(a)/h\right)\right]
=+\infty.
\end{equation}
\end{enumerate}
\end{theorem}

The critical case $\|a\|_\infty=1/2$ is not covered by the theorem.

We prove the asymptotic identity \eqref{eq:nonlinear-limit} by establishing
its lower bound in Section~\ref{sec:lower} and the matching upper bound in
Section~\ref{sec:upper}. Here we prove only the supercritical assertion
\eqref{eq:nonlinear-supercritical}.

\begin{proof}[Proof of the supercritical assertion]
Choose $(t_0,x_0)$ with $a(t_0,x_0)>1/2$, and a closed interval
$I\subset(0,1)$ on which $a(t,x_0)>1/2$. Put
$a_*:=\min_{t\in I}a(t,x_0)>1/2$, and choose $z>0$ such that
$2a_*\sinc z>1$. Next, choose a nonzero $\eta\in C_c^\infty(I)$
and a smooth function $m$ with $m(0)=0$ and $m=x_0$ on a
neighborhood of $I$. Fix $\varepsilon>0$, to be chosen sufficiently
small, and for $h>0$ set
\[
 g_h(t):=m(t)+\varepsilon\eta(t)\sin(zt/h),\qquad
 v_h(t):=g_h'(t),\quad t\in[0,1].
\]
Define the probability measure $\mathbb Q_h$ by
\[
 \frac{d\mathbb Q_h}{d\Pp}
 :=\exp\left(\int_0^1v_h(t)\,dW_t
           -\frac12\int_0^1v_h^2(t)\,dt\right).
\]
By Girsanov's theorem, $B_t:=W_t-g_h(t)$ is a $\mathbb Q_h$
Brownian motion. 

Next, set $E_h:=\{\sup_{0\le t\le1}|B_t|<\varepsilon\}$. Its
$\mathbb Q_h$-probability is
$p_\varepsilon:=\Pp(\sup_{0\le t\le1}|W_t|<\varepsilon)>0$,
independent of $h$. Introduce the probability measure $\mathbb R_h$ by
\[
 \frac{d\mathbb R_h}{d\Pp}
 :=\frac{\one_{E_h}}{p_\varepsilon}
 \exp\left(\int_0^1v_h(t)\,dW_t
           -\frac12\int_0^1v^2_h(t)\,dt\right)
\]
where $\one_{\cdot}$ denotes the indicator function. 

The relative entropy $H(\mathbb R_h\mid\Pp):=
\E_{\mathbb R_h}[\log(d\mathbb R_h/d\Pp)]$ is given by
\[
 H(\mathbb R_h\mid\Pp)
 =-\log p_\varepsilon+\frac12\int_0^1v_h(t)^2\,dt
    +\E_{\mathbb R_h}\left[\int_0^1v_h(t)\,dB_t\right].
\]
The last expectation is zero: under $\mathbb Q_h$, both the Brownian
law and $E_h$ are invariant under $B\mapsto-B$, whereas
$\int_0^1v_h(t)\,dB_t$ changes sign. Thus
\[
 H(\mathbb R_h\mid\Pp)
 =\frac12\int_0^1\left(m'(t)+\varepsilon\eta'(t)\sin(zt/h)
       +\frac{\varepsilon z}{h}\eta(t)\cos(zt/h)\right)^2dt
       -\log p_\varepsilon.
\]
Since $m'\equiv0$ on $I$, the Riemann--Lebesgue lemma yields
\begin{equation}\label{eq:super-entropy}
 \lim_{h\downarrow0}h^2H(\mathbb R_h\mid\Pp)
 =\frac{\varepsilon^2z^2}{4}\int_I\eta(t)^2\,dt.
\end{equation}

We next estimate $\E_{\mathbb R_h}\left[N_h(a)\right]$. Under
$\mathbb Q_h$, write $dW_t=dB_t+v_h(t)\,dt$ and put
\[
 A_t:=\int_{\tau_h(t)}^t a(s,W_s)\,dB_s,\qquad
 D_t:=\int_{\tau_h(t)}^t a(s,W_s)v_h(s)\,ds.
\]
Then
\begin{equation}\label{eq:super-decomp}
 N_h(a)=\int_0^1 A_t\,dB_t+\int_0^1D_t\,dB_t
          +\int_0^1A_tv_h(t)\,dt+\int_0^1D_tv_h(t)\,dt.
\end{equation}
As $a$ is globally Lipschitz and $\sup_{h,t}|g_h(t)|<\infty$,
there is a constant $C<\infty$, independent of $h\in(0,1)$,
such that
\[
 \sup_{h\in(0,1)}\sup_{s\in[0,1]}
 \E_{\mathbb Q_h}\left[a(s,W_s)^2\right]\le C.
\]
Moreover, $\|v_h\|_\infty\le C/h$. By It\^o's isometry and
Cauchy--Schwarz,
\[
 \E_{\mathbb Q_h}\left[\int_0^1 A_t^2\,dt\right]\le Ch,\qquad
 \E_{\mathbb Q_h}\left[\int_0^1 D_t^2\,dt\right]\le C,\qquad
 \int_0^1v_h(t)^2\,dt\le Ch^{-2}.
\]
Consequently, the expected absolute values under $\mathbb Q_h$
of the first three terms in \eqref{eq:super-decomp} are bounded,
respectively, by $C\sqrt h$, $C$, and $C/\sqrt h$. Since
$d\mathbb R_h/d\mathbb Q_h=\one_{E_h}/p_\varepsilon$,
the same bounds under $\mathbb R_h$ hold with an additional
factor $p_\varepsilon^{-1}$. In particular,
\begin{equation}\label{eq:super-cross}
 \lim_{h\downarrow0}h\left(
 \E_{\mathbb R_h}\left[N_h(a)\right]
 -\E_{\mathbb R_h}\left[\int_0^1D_tv_h(t)\,dt\right]
 \right)=0.
\end{equation}

Put
\[
 q_h(t):=v_h(t)-m'(t)
 =\varepsilon\eta'(t)\sin(zt/h)
   +\frac{\varepsilon z}{h}\eta(t)\cos(zt/h).
\]
Since $q_h$ is supported in $I$ and $m'$ vanishes on a
neighborhood of $I$, the supports of $q_h$ and $m'$ are separated
by a positive distance. The mixed terms involving $m'$ and $q_h$
in the double integral therefore vanish whenever $h$ is smaller
than that distance. On $E_h$, the term involving only $m'$ is
$O(h)$, uniformly in $h$: indeed $W=B+g_h$ is uniformly bounded
there, $a$ is Lipschitz, and the integration strip has area at
most $h$.
It follows that
\begin{equation}\label{eq:super-drift-reduction}
 \lim_{h\downarrow0}h\left(
 \E_{\mathbb R_h}\left[\int_0^1D_tv_h(t)\,dt\right]
 -\E_{\mathbb R_h}\left[
 \int_0^1\int_{\tau_h(t)}^t a(s,W_s)q_h(s)q_h(t)\,ds\,dt
 \right]\right)=0.
\end{equation}

Let $L$ be a global Lipschitz constant of $a$. On $E_h$,
for $s$ in the support of $q_h$ we have $m(s)=x_0$ and hence
\[
 |W_s-x_0|\le\varepsilon(1+\|\eta\|_\infty),\qquad
 |a(s,W_s)-a(s,x_0)|
 \le L\varepsilon(1+\|\eta\|_\infty).
\]
The Cauchy--Schwarz inequality implies 
\[
 \int_0^1\int_{\tau_h(t)}^t
 |q_h(s)q_h(t)|\,ds\,dt
 \le h\int_0^1q_h(t)^2\,dt.
\]
Thus
\begin{align}
 &h\E_{\mathbb R_h}\left[
 \int_0^1\int_{\tau_h(t)}^t a(s,W_s)q_h(s)q_h(t)\,ds\,dt
 \right]\notag\\
 &\quad\ge
 h\int_0^1\int_{\tau_h(t)}^t a(s,x_0)q_h(s)q_h(t)\,ds\,dt
 -L\varepsilon(1+\|\eta\|_\infty)
 h^2\int_0^1q_h(t)^2\,dt.\label{eq:super-freeze}
\end{align}
Directly from the definition of $q_h$,
\begin{equation}\label{eq:super-qnorm}
 \lim_{h\downarrow0}h^2\int_0^1q_h(t)^2\,dt
 =\frac{\varepsilon^2z^2}{2}\int_I\eta(t)^2\,dt.
\end{equation}
Also, since
$q_h(t)=\varepsilon z\eta(t)\cos(zt/h)/h+O(1)$ uniformly in $t$
(for fixed $\varepsilon$), the lower-order terms contribute $o(1)$
after multiplication of the double integral by $h$. Setting
$t=s+hu$ gives
\begin{align*}
 &\lim_{h\downarrow0}h
 \int_0^1\int_{\tau_h(t)}^t
 a(s,x_0)q_h(s)q_h(t)\,ds\,dt\\
 &\quad=\varepsilon^2z^2\lim_{h\downarrow0}
 \int_0^1\!\int_0^1
 a(s,x_0)\eta(s)\eta(s+hu)
 \cos(zs/h)\cos(zs/h+zu)\,du\,ds\\
 &\quad=\frac{\varepsilon^2z^2}{2}\sinc z
 \int_Ia(s,x_0)\eta(s)^2\,ds.
\end{align*}
Here $\eta$ is extended by zero outside $[0,1]$; the
indicator $s+hu\le1$ may be omitted for all sufficiently small
$h$ on the support of $\eta$. Indeed,
$\eta(s+hu)\to\eta(s)$ uniformly in $s,u\in[0,1]$.
The last equality follows from
$\cos x\cos(x+zu)=\frac12\cos(zu)+\frac12\cos(2x+zu)$,
the Riemann--Lebesgue lemma applied to the remaining fixed
$L^1$ amplitude, and $\int_0^1\cos(zu)\,du=\sinc z$.
Combining this limit with \eqref{eq:super-cross},
\eqref{eq:super-drift-reduction}, \eqref{eq:super-freeze}, and
\eqref{eq:super-qnorm} yields
\begin{equation}\label{eq:super-energy}
 \liminf_{h\downarrow0}h\E_{\mathbb R_h}\left[N_h(a)\right]
 \ge \frac{\varepsilon^2z^2}{2}\int_I
 \bigl(a(t,x_0)\sinc z-L\varepsilon(1+\|\eta\|_\infty)\bigr)
 \eta(t)^2\,dt.
\end{equation}

Finally, the entropy inequality gives
\[
 h\log\E\left[e^{N_h(a)/h}\right]
 \ge\E_{\mathbb R_h}\left[N_h(a)\right]
       -hH(\mathbb R_h\mid\Pp).
\]
By \eqref{eq:super-entropy} and \eqref{eq:super-energy},
\[
 \liminf_{h\downarrow0}h^2\log\E\left[e^{N_h(a)/h}\right]
 \ge\frac{\varepsilon^2z^2}{4}\int_I
 \bigl(2a(t,x_0)\sinc z-1
       -2L\varepsilon(1+\|\eta\|_\infty)\bigr)\eta(t)^2\,dt.
\]
Since $2a_*\sinc z>1$ and $\eta$ is nonzero, the right side
is strictly positive when $\varepsilon>0$ is sufficiently small.
This proves \eqref{eq:nonlinear-supercritical}.
\end{proof}

We conclude this section with the following corollary, which describes the asymptotic
performance of the delayed-delta hedging strategy.
\begin{corollary}
Assume that $f\in C^4(\R)$ is convex and that its third and fourth derivatives
are uniformly bounded. If $\sup_{x\in\R} f''(x)<1/2$, then
\begin{equation}\label{eq:main}
\lim_{h\downarrow0}h\log\E\left[e^{\frac{1}{h}\left(f(W_1)-\int_{0}^1 u_x\left(\tau_h(t), W_{\tau_h(t)}\right)dW_t\right)}\right]
= u(0,0)+\int_0^1\sup_{x\in\R}\Psi\bigl(u_{xx}(t,x)\bigr)\,dt.
\end{equation}
If $\sup_{x\in\R} f''(x)>1/2$, then
\begin{equation}\label{eq:main-supercritical}
\lim_{h\downarrow0}h\log\E\left[e^{\frac{1}{h}\left(f(W_1)-\int_{0}^1 u_x\left(\tau_h(t), W_{\tau_h(t)}\right)dW_t\right)}\right]
=+\infty.
\end{equation}
\end{corollary}
\begin{proof}
Let $a(t,x):=u_{xx}(t,x)=\E\left[f''(x+W_{1-t})\right]$. 
Here the heat-semigroup representation extends to $t=1$, with $a(1,x)=f''(x)$.
Since $f$ is convex, $a\ge0$.
Moreover, the boundedness of $f^{(3)}$, $f^{(4)}$ and the identity $a_t=-\frac{1}{2}a_{xx}$ imply that
\[
|a(t,x)-a(s,y)|\le
\|f^{(3)}\|_\infty|x-y|+
\frac12\|f^{(4)}\|_\infty|t-s|,
\]
so $a$ is globally Lipschitz. By \eqref{eq:ito-identity},
\[
f(W_1)-\int_0^1u_x\bigl(\tau_h(t),W_{\tau_h(t)}\bigr)\,dW_t
=u(0,0)+N_h(a).
\]
The heat-semigroup representation and $a(1,x)=f''(x)$ give
$\|a\|_\infty=\sup_{x\in\R}f''(x)$, including the value $+\infty$.
The assertions \eqref{eq:main} and \eqref{eq:main-supercritical} follow from
\eqref{eq:nonlinear-limit} and \eqref{eq:nonlinear-supercritical}, respectively.
\end{proof}

\section{Auxiliary lemmas}\label{sec:auxiliary}
We collect the Gaussian comparison and exponential moment estimates
used in the lower- and upper-bound proofs of Theorem~\ref{thm:nonlinear}.

\begin{lemma}\label{lem:domination}
Let $I:=[r_0,r_1]\subset[0,1]$ and let $k\in L^2([0,1]\times I)$ be deterministic.
Let $\beta$ be predictable with respect to the Brownian filtration, and suppose
that $|\beta_s|\le M$ for $ds\otimes d\Pp$-almost every $(s,\omega)$,
where $M\ge0$ is a finite constant. For almost every $t\in I$, set
\[
Z_t:=\int_0^1 k(s,t)\beta_s\,dW_s,\qquad
\widetilde Z_t:=M\int_0^1 k(s,t)\,dW_s.
\]
Then
\begin{equation}\label{eq:gaussiancomp}
\E\left[\exp\left(\frac12\int_I Z_t^2\,dt\right)\right]
\le \E\left[\exp\left(\frac12\int_I \widetilde Z_t^2\,dt\right)\right],
\end{equation}
where the right-hand side may be infinite.
\end{lemma}
\begin{proof}
Let
$(e_j)_{j\ge1}$ be an orthonormal basis of $L^2(I)$. By It\^o's
isometry and Tonelli's theorem, $Z$ and $\widetilde Z$ belong to $L^2(I)$
almost surely. For each $j$ put
\[
w_j(s):=\int_I k(s,t)e_j(t)\,dt.
\]
Then $w_j\in L^2([0,1])$, and stochastic Fubini gives
\[
X_j:=\langle Z,e_j\rangle_{L^2(I)}
 =\int_0^1 w_j(s)\beta_s\,dW_s,\qquad
Y_j:=\langle\widetilde Z,e_j\rangle_{L^2(I)}
 =M\int_0^1 w_j(s)\,dW_s.
\]

Fix $d\ge1$. For $g\in\R^d$, define $f_g:=\sum_{j=1}^d g_jw_j$ and
$N_u:=\int_0^u f_g(s)\beta_s\,dW_s$, $u\in [0,1]$. Its quadratic variation satisfies
$\langle N\rangle_1\le M^2\int_0^1 f_g(s)^2\,ds$. The stochastic
exponential $\exp(N_u-\langle N\rangle_u/2)$ is a nonnegative local
martingale starting at one, so its expectation at time one is at most
one. It follows that
\[
\E\left[e^{\sum_{j=1}^d g_jX_j}\right]
=\E\left[e^{N_1}\right]
\le \exp\left(\frac{M^2}{2}\int_0^1 f_g(s)^2\,ds\right)
=\E\left[e^{\sum_{j=1}^d g_jY_j}\right],
\]
where the last equality uses the Gaussian moment-generating function
of $M\int_0^1 f_g(s)\,dW_s$.

Let $\gamma_d$ be the standard Gaussian measure on $\R^d$. Integrating
the preceding inequality against $\gamma_d$ and applying Tonelli's theorem
and the identity $\int_{\R^d}e^{g\cdot x}\,\gamma_d(dg)=e^{|x|^2/2}$
gives
\[
\E\left[e^{\frac12\sum_{j=1}^d X_j^2}\right]
\le \E\left[e^{\frac12\sum_{j=1}^d Y_j^2}\right].
\]
Parseval's identity and monotone convergence yield
\eqref{eq:gaussiancomp} as $d\to\infty$.
\end{proof}

\begin{lemma}\label{lem:moment}
Fix $h>0$ and $I:=[r,r+\Delta]\subset[0,1]$ with $\Delta>0$.
Let $k$ be a deterministic kernel on $I^2$ supported on $0<t-s<h$ and
let $\beta$ be predictable, with essential bounds $\|k\|_\infty\le c$
and $\|\beta\|_\infty\le M$. Set
\[
Z_t:=\int_r^t k(s,t)\beta_s\,dW_s,\qquad J_h:=\int_I Z_t\,dW_t.
\]
If $2Mc<1$, then
\begin{equation}\label{eq:moment}
h\log\E\left[e^{J_h/h}\mid\F_r\right]
\le \frac{\Delta M^2c^2}{1-4M^2c^2}
\qquad\text{almost surely}.
\end{equation}
\end{lemma}
\begin{proof}
Put $\mathcal V:=\int_I Z_t^2\,dt$. The stochastic exponential of the martingale $(\frac{2}{h}\int_r^u Z_t\,dW_t)_{r\le u\le r+\Delta}$ is a
nonnegative supermartingale starting at one at time $r$. Its terminal
value is $D:=\exp(\frac{2J_h}{h}-\frac{2\mathcal V}{h^2})$. By conditional Cauchy--Schwarz,
\begin{equation}\label{eq:cs}
\E\left[e^{J_h/h}\mid\F_r\right]
=\E\left[D^{1/2}e^{\mathcal V/h^2}\mid\F_r\right]
\le \E\left[e^{2\mathcal V/h^2}\mid\F_r\right]^{1/2}.
\end{equation}

Put $\widetilde Z_t:=M\int_I k(s,t)\,dW_s$ and
$\widetilde{\mathcal V}:=\int_I\widetilde Z_t^2\,dt$. For any $A\in\F_r$ of positive
probability, under $\Pp(\cdot\mid A)$ the process
$(W_{r+u}-W_r)_{0\le u\le\Delta}$ is Brownian with respect to
$(\F_{r+u})_{0\le u\le\Delta}$ by independence of increments.
The proof of Lemma~\ref{lem:domination} applies to any filtration with
respect to which the integrator is Brownian. Apply it under
$\Pp(\cdot\mid A)$ on the shifted interval with kernel $2k/h$.
The law of $\widetilde Z$ is unchanged under this conditioning.
Since $A$ is arbitrary, we obtain
\begin{equation}\label{eq:moment-conditional-comp}
\E\left[e^{2\mathcal V/h^2}\mid\F_r\right]\le\E\left[e^{2\widetilde{\mathcal V}/h^2}\right]
\qquad\text{almost surely}.
\end{equation}

Let $V$ be the integral operator on $L^2(I)$ with
$(Vv)(t):=\int_I k(s,t)v(s)\,ds$. The support has length at most $h$ in
each variable with the other fixed. Schur's test and integration give
\begin{equation*}
\|V\|_{\mathrm{op}}\le ch,\qquad \|V\|_{\HS}^2\le c^2\Delta h.
\end{equation*}
By It\^o's isometry, $\widetilde Z$ has covariance $M^2VV^*$, whose
eigenvalues $\mu_j$ satisfy
\begin{equation*}
0\le\mu_j\le M^2c^2h^2,\qquad \sum_j\mu_j\le M^2c^2\Delta h.
\end{equation*}
Set $\rho:=4M^2c^2<1$. In an eigenbasis, the Gaussian coordinates are
independent, so their squared exponential moments and monotone
convergence give
\begin{equation*}
\E\left[e^{2\widetilde{\mathcal V}/h^2}\right]=\prod_j(1-4\mu_j/h^2)^{-1/2}.
\end{equation*}
Using $-\log(1-x)\le x/(1-\rho)$ for $0\le x\le\rho$ in this product,
we obtain
\[
\log\E\left[e^{2\widetilde{\mathcal V}/h^2}\right]
\le \frac{2}{h^2(1-\rho)}\sum_j\mu_j
\le \frac{2\Delta M^2c^2}{h(1-4M^2c^2)}.
\]
Equations \eqref{eq:cs} and \eqref{eq:moment-conditional-comp} now give
the upper bound in \eqref{eq:moment}.
\end{proof}

\section{Lower bound}\label{sec:lower}
We start with the following lemma.

\begin{lemma}[Weak convergence under the tilted measure]\label{lem:tube}
Fix a deterministic continuous function $b:[0,1]\to[0,1/2)$. By Proposition~\ref{prop:deterministic}, the exponential moment in the denominator below is finite, so
\begin{equation}\label{eq:tilt}
\frac{d\mathbb Q_h^{b}}{d\Pp}
:=\frac{\exp\{N_h(b)/h\}}{\E\left[\exp\{N_h(b)/h\}\right]}
\end{equation}
defines a probability measure.
Under $\mathbb Q_h^{b}$, the coordinate process $W$ converges in distribution in $C_0([0,1])$, equipped with the uniform norm, to
\[
Y_t:=\int_0^t(1-2b(s))^{-1/2}\,dB_s,\qquad 0\le t\le1,
\]
where $B$ is a standard Brownian motion under $\Pp$ (for example, $B=W$). Here $C_0([0,1])$ denotes the space of continuous paths starting at zero. In particular, for every $\varphi\in C_0([0,1])$ and $\rho>0$,
\begin{equation}\label{eq:tube}
\lim_{h\downarrow0}h\log\E\left[e^{N_h(b)/h}
\one_{\{\|W-\varphi\|_\infty<\rho\}}\right]
=\int_0^1\Psi(b(t))\,dt.
\end{equation}
\end{lemma}
\begin{proof}
The key observation is that the exponential change of measure changes the
covariance of the Brownian motion but preserves its Gaussian character.
To see this, let $(e_{h,j})_{j\ge1}$ be an orthonormal eigenbasis of
$A_{h,b}$ (including an orthonormal basis of its kernel), with corresponding
 eigenvalues $\lambda_{h,j}$, and write
$\xi_{h,j}=\int_0^1 e_{h,j}(s)\,dW_s$. Under $\Pp$, these variables are
independent standard Gaussians. By the spectral expansion of the kernel,
It\^o's formula and It\^o's isometry,
\[
S_n:=\frac12\sum_{j=1}^n\lambda_{h,j}(\xi_{h,j}^2-1)
\longrightarrow \frac{N_h(b)}h\qquad\text{in }L^2(\Pp).
\]
Recall from \eqref{eq:det-schur} that
$|\lambda_{h,j}|\le 2\|b\|_\infty<1$. We may therefore choose $p>1$
so that $2p\|b\|_\infty<1$. The Gaussian moment formula gives
\[
\log\E[e^{pS_n}]
=\frac12\sum_{j=1}^n\bigl[-\log(1-p\lambda_{h,j})-p\lambda_{h,j}\bigr]
\le C\sum_{j\ge1}\lambda_{h,j}^2<\infty,
\]
where $C$ depends only on $p$ and $\|b\|_\infty$. Thus $e^{S_n}$ is
uniformly integrable and converges in $L^1(\Pp)$ to $e^{N_h(b)/h}$.

For the finite sum $S_n$, the effect of the change of measure is
particularly simple. Multiplying the standard Gaussian density of
$\xi_{h,j}$ by the normalized exponential factor
$(1-\lambda_{h,j})^{1/2}\exp(\lambda_{h,j}\xi_{h,j}^2/2)$ replaces its
variance by $(1-\lambda_{h,j})^{-1}$ and leaves its mean unchanged.
Since the density factors over the coordinates, independence is preserved.
The $L^1$ convergence established above allows us to pass to the limit
in the normalized densities. Consequently, under $\mathbb Q_h^b$ the
variables $(\xi_{h,j})_{j\ge1}$ are independent centered Gaussians with
variances $(1-\lambda_{h,j})^{-1}$.

For any $g\in L^2([0,1])$, expand the Wiener integral in this basis:
\[
\int_0^1 g(u)\,dW_u
=\sum_{j\ge1}\langle g,e_{h,j}\rangle_{L^2}\xi_{h,j}.
\]
The series converges in $L^2(\Pp)$ by It\^o's isometry and in
$L^2(\mathbb Q_h^b)$ because the tilted variances are uniformly bounded
by $(1-2\|b\|_\infty)^{-1}$. The two limits agree by absolute continuity.
Taking $g$ to be linear combinations of indicators of intervals $[0,t]$,
we conclude that $W$ is a centered Gaussian process under $\mathbb Q_h^b$,
with covariance
\begin{equation}\label{eq:tiltcov}
\Gamma_h(s,t)=\left\langle\one_{[0,s]},
(I-A_{h,b})^{-1}\one_{[0,t]}\right\rangle_{L^2}.
\end{equation}

It remains to understand this covariance as $h\downarrow0$. For a
continuous function $g$ and $t\in(0,1)$, continuity of $b$ and $g$ gives
\[
(A_{h,b}g)(t)=\frac1h\int_{(t-h)\vee0}^{(t+h)\wedge1}
b(s\wedge t)g(s)\,ds\longrightarrow 2b(t)g(t).
\]
The integrals are bounded in absolute value by
$2\|b\|_\infty\|g\|_\infty$, so the convergence also holds in $L^2$.
Let $M_b$ denote multiplication by $2b$. By density of continuous
functions in $L^2$ and the uniform operator norm bound
$\|A_{h,b}\|_{\mathrm{op}},\|M_b\|_{\mathrm{op}}\le2\|b\|_\infty<1$,
we have $A_{h,b}\to M_b$ strongly. The resolvent identity
\[
(I-A_{h,b})^{-1}-(I-M_b)^{-1}
=(I-A_{h,b})^{-1}(A_{h,b}-M_b)(I-M_b)^{-1}
\]
then gives strong convergence of the inverses, since their norms are
bounded by $(1-2\|b\|_\infty)^{-1}$. In particular,
\[
\Gamma_h(s,t)\longrightarrow\int_0^{s\wedge t}\frac{du}{1-2b(u)}.
\]
This is precisely the covariance of $Y$. As the processes are centered
Gaussian, we obtain convergence of the finite-dimensional distributions.
To upgrade this to weak convergence on the path space, note that
\eqref{eq:tiltcov} and the bound on the inverse imply
\[
\E_{\mathbb Q_h^b}[(W_t-W_s)^2]
\le\frac{|t-s|}{1-2\|b\|_\infty}.
\]
By Gaussianity,
\[
\E_{\mathbb Q_h^b}[(W_t-W_s)^4]
\le\frac{3|t-s|^2}{(1-2\|b\|_\infty)^2}.
\]
Together with $W_0=0$, this gives tightness by Kolmogorov's criterion,
and hence weak convergence to $Y$ in $C_0([0,1])$.

Finally, $Y$ has the law of Brownian motion run at the deterministic
clock $V(t)=\int_0^t(1-2b(u))^{-1}\,du$. Since $V$ is continuous and
strictly increasing, the full support of Brownian motion implies that
$Y$ has full support in $C_0([0,1])$: indeed, one can apply the Brownian
support property to the path $\varphi\circ V^{-1}$ on $[0,V(1)]$.
The tube $\{w:\|w-\varphi\|_\infty<\rho\}$ is open, so weak convergence
and the Portmanteau theorem give
\[
\liminf_{h\downarrow0}\mathbb Q_h^b(\|W-\varphi\|_\infty<\rho)
\ge\Pp(\|Y-\varphi\|_\infty<\rho)>0.
\]
Thus restricting to the tube has no cost on the exponential scale:
$h\log\mathbb Q_h^b(\|W-\varphi\|_\infty<\rho)\to0$.
Returning to the original measure via \eqref{eq:tilt}, we obtain
\begin{align*}
h\log\E\left[e^{N_h(b)/h}\one_{\{\|W-\varphi\|_\infty<\rho\}}\right]
={}&h\log\E[e^{N_h(b)/h}]\\
&+h\log\mathbb Q_h^b(\|W-\varphi\|_\infty<\rho).
\end{align*}
The first term converges to $\int_0^1\Psi(b(t))\,dt$ by
Proposition~\ref{prop:deterministic}, while the second tends to zero.
This proves \eqref{eq:tube}.
\end{proof}

We conclude this section by establishing the lower bound in \eqref{eq:nonlinear-limit}.
 \begin{proof}
Throughout the rest of this section, assume
$0\le a(t,x)\le\kappa<1/2$ and let $L$ denote its global
Lipschitz constant. We prove
\begin{equation*}
\liminf_{h\downarrow0}h\log\E\left[e^{N_h(a)/h}\right]\ge \int_0^1\sup_{x\in\R}\Psi(a(t,x))\,dt.
\end{equation*}
Fix a continuous deterministic path $\varphi$ with $\varphi(0)=0$ and a radius $\rho>0$. Define
\[
b(t):=a(t,\varphi(t)),\qquad
A_\rho:=\{\|W-\varphi\|_\infty<\rho\}.
\]
Then $b$ is continuous and $0\le b\le\kappa$. Let
\[
\tau:=\inf\{t\in[0,1]:|W_t-\varphi(t)|\ge\rho\},
\]
with $\tau:=\infty$ if there is no exit. Define
\[
\beta_s:=\big(a(s,W_s)-b(s)\big)\one_{\{s\le\tau\}},\qquad
\mathcal E_h:=\int_0^1\int_{\tau_h(t)}^t\beta_s\,dW_s\,dW_t.
\]
The process $\beta$ is predictable and satisfies $|\beta_s|\le L\rho$.
To see why the stopping does not affect the error on $A_\rho$, note that
by the stopping property of It\^o integrals,
\[
\int_0^{t\wedge\tau}(a(s,W_s)-b(s))\,dW_s
=\int_0^{t\wedge\tau}\beta_s\,dW_s,
\qquad t\in[0,1],
\]
almost surely, simultaneously for all $t$. Taking differences at $t$
and $\tau_h(t)$ shows that the corresponding inner integrals agree
for $t\le\tau$. Applying the stopping property once more, now to
the outer integral, gives
\begin{align*}
&\int_0^{1\wedge\tau}\int_{\tau_h(t)}^t
(a(s,W_s)-b(s))\,dW_s\,dW_t\\
&\qquad=\int_0^{1\wedge\tau}\int_{\tau_h(t)}^t
\beta_s\,dW_s\,dW_t.
\end{align*}
Since $\tau=\infty$ on $A_\rho$, linearity of stochastic integration
then yields, almost surely,
\begin{equation*}
N_h(a)=N_h(b)+\mathcal E_h\qquad\text{on }A_\rho.
\end{equation*}

Fix $\alpha\in(0,1)$ and choose $\rho$ sufficiently small that
$2\alpha L\rho/(1-\alpha)<1$. On $A_\rho$, write
\[
e^{\alpha N_h(b)/h}
=e^{\alpha N_h(a)/h}e^{-\alpha \mathcal E_h/h}.
\]
H\"older's inequality with exponents $1/\alpha$ and $1/(1-\alpha)$ yields
\begin{align*}
\E\left[e^{\alpha N_h(b)/h}\one_{A_\rho}\right]
&\le\bigl(\E\left[e^{N_h(a)/h}\right]\bigr)^\alpha
\left(\E\left[e^{-\alpha \mathcal E_h/((1-\alpha)h)}\right]\right)^{1-\alpha}.
\end{align*}
Taking logarithms and rearranging (with the inequality automatic if the moment of $N_h(a)$ is infinite),
\begin{align*}
h\log\E\left[e^{N_h(a)/h}\right]
&\ge\frac h\alpha\log\E\left[e^{\alpha N_h(b)/h}\one_{A_\rho}\right]\\
&\quad-\frac{1-\alpha}{\alpha}h\log\E\left[e^{-\alpha \mathcal E_h/((1-\alpha)h)}\right].
\end{align*}
Since $\|\alpha b\|_\infty<1/2$, we may apply Lemma~\ref{lem:tube} to $\alpha b$ in the first logarithm. Applying Lemma~\ref{lem:moment} to the second gives
\begin{align*}
\liminf_{h\downarrow0}h\log\E\left[e^{N_h(a)/h}\right]
\ge{}&\frac1\alpha\int_0^1\Psi(\alpha b(t))\,dt\\
&-\frac{1-\alpha}{\alpha}
\frac{(\alpha L\rho/(1-\alpha))^2}{1-4(\alpha L\rho/(1-\alpha))^2}.
\end{align*}
The left-hand side is independent of $\rho$ and $\alpha$. With $\alpha$ fixed, let $\rho\downarrow0$. The second term tends to zero. Next let $\alpha\uparrow1$. Continuity of $\Psi$ on $[0,\kappa]$ and bounded convergence imply
\begin{equation*}
\liminf_{h\downarrow0}h\log\E\left[e^{N_h(a)/h}\right]
\ge\int_0^1\Psi(a(t,\varphi(t)))\,dt.
\end{equation*}

This, together with the equality
\[
\sup_{\varphi\in C_0([0,1])}\int_0^1\Psi(a(t,\varphi(t)))\,dt
=\int_0^1\sup_{x\in\R}\Psi(a(t,x))\,dt,
\]
gives the desired lower bound. Indeed, uniform continuity in $t$, uniformly in $x$,
allows the right-hand side to be approximated by choosing a nearly
maximizing constant on each interval of a sufficiently fine partition.
Connecting these constants linearly on intervals of arbitrarily small
total length, including a transition from zero at the origin, produces
continuous paths with arbitrarily small loss, since $\Psi$ is bounded
on $[0,\kappa]$. 
\end{proof}

\section{Upper bound by block estimates}\label{sec:upper}
Assume $0\le a\le\kappa<1/2$, and let $L$ be a Lipschitz constant of $a$.
We divide $[0,1]$ into intervals of length much larger than $h$ but
small enough that $a(s,W_s)$ can be approximated by its value at the
left endpoint. We first deal with the terms in $N_h(a)$ that cross
interval boundaries. Fix $r>1$ such that $2r^2\kappa<1$, and set
\[
q=\frac{r}{r-1},\qquad \delta=h^{3/4},\qquad \rho=h^{1/8}.
\]
We keep $r$ fixed until the end of the proof. All error estimates below
are deterministic and uniform in the partition shift, the interval and
the Brownian past.

For $v\in[0,\delta)$, partition $[0,1]$ at the points $v+j\delta$ lying in $(0,1)$, with $j\in\mathbb Z$.
For each interval $I=[u,w]$ in the partition, put
\[
J_I=\int_u^w\int_{\max\{u,t-h\}}^t a(s,W_s)\,dW_s\,dW_t,
\qquad N_h^v=\sum_I J_I.
\]
The reason for introducing the shift $v$ is that, for any pair
$0<t-s<h<\delta$, the set of shifts that separate $s$ and $t$ has
length $t-s$. Averaging over all shifts and applying stochastic Fubini
(justified by It\^o's isometry and boundedness of $a$), we get
\begin{equation*}
N_h(a)=\frac1\delta\int_0^\delta N_h^v\,dv+D_h,\qquad
D_h=\int_0^1\int_{(t-h)^+}^t\frac{t-s}{\delta}a(s,W_s)\,dW_s\,dW_t.
\end{equation*}
By Lemma~\ref{lem:moment},
$h\log\E[e^{qD_h/h}]=o(1)$, since the coefficient is bounded by $q\kappa h/\delta\to0$.
H\"older's and Jensen's inequalities therefore give
\begin{equation}\label{eq:upper-error-moment}
h\log\E[e^{N_h(a)/h}]
\le \frac1r\sup_{0\le v<\delta}h\log\E[e^{rN_h^v/h}]+o(1).
\end{equation}
It therefore suffices to estimate the exponential moments on the
individual intervals, conditional on the Brownian past.

Let $I=[u,u+\delta]$ be an interval of length $\delta$, let $a_u=a(u,W_u)$, and define
\[
J_I^{(1)}:=\int_u^{u+\delta}(W_t-W_{(t-h)\vee u})\,dW_t.
\]
By rescaling $I$ to $[0,1]$ and applying Proposition~\ref{prop:deterministic} with constant coefficient $c$, we obtain
\begin{equation}\label{eq:upper-uniform}
\frac h\delta\log\E[e^{cJ_I^{(1)}/h}\mid\F_u]
=\Psi(c)+o(1),\qquad 0\le c\le r^2\kappa.
\end{equation}
Here $h/\delta\to0$. The convergence is uniform in $c$: the functions
on the left are nondecreasing in $c\ge0$ by convexity and the fact
that the integral has mean zero, and their limit is continuous.

We can replace $a(s,W_s)$ by $a_u$ as long as Brownian motion stays
close to $W_u$. Let $\sigma$ be the first time in $I$ at which
$|W_s-W_u|\ge\rho$, with $\sigma=u+\delta$ if there is no such time.
Write $\mathcal E_I$ for the double integral defining $J_I$ with coefficient
$(a(s,W_s)-a_u)\one_{\{s\le\sigma\}}$ in place of $a(s,W_s)$.
By the Lipschitz property of $a$ and the definition of $\sigma$, this
predictable coefficient is bounded by $L(\delta+\rho)$.
On the event $A_I=\{\sup_{s\in I}|W_s-W_u|<\rho\}$, locality of stochastic
integration gives $J_I=a_uJ_I^{(1)}+\mathcal E_I$. Applying H\"older's inequality conditionally on $\F_u$ and using
Lemma~\ref{lem:moment} together with \eqref{eq:upper-uniform}, we obtain
\begin{align*}
\E[e^{rJ_I/h}\one_{A_I}\mid\F_u]
&\le \E[e^{r^2a_uJ_I^{(1)}/h}\mid\F_u]^{1/r}
       \E[e^{qr\mathcal E_I/h}\mid\F_u]^{1/q}\\
&\le \exp\left\{\frac\delta h
\left[\frac1r\Psi(r^2a_u)+o(1)\right]\right\}.
\end{align*}
The error is uniform because $L(\delta+\rho)\to0$ and
\eqref{eq:upper-uniform} holds uniformly in $c$.
On the complementary event, Brownian motion moves by at least $\rho$
during $I$. Its probability is small enough to compensate for the
exponential moment of $J_I$. By conditional H\"older's inequality,
Lemma~\ref{lem:moment} and the reflection principle,
\[
\E[e^{rJ_I/h}\one_{A_I^c}\mid\F_u]
\le \exp\left\{\frac\delta{rh}\frac{r^4\kappa^2}{1-4r^4\kappa^2}\right\}
\left(4e^{-\rho^2/(2\delta)}\right)^{1/q}=o(1).
\]
The right-hand side tends to zero because $\rho^2h/\delta^2=h^{-1/4}\to\infty$.
Adding the estimates on $A_I$ and $A_I^c$, we obtain
\begin{equation}\label{eq:upper-cell}
h\log\E[e^{rJ_I/h}\mid\F_u]
\le \frac\delta r\sup_{x\in\R}\Psi(r^2a(u,x))+o(\delta).
\end{equation}
Indeed, $\Psi\ge0$ and the contribution of $A_I^c$ tends to zero;
taking the scaled logarithm of the sum adds at most
$h\log 2=o(\delta)$ to the larger of the two scaled logarithms. For the at most two shorter intervals at the ends
of the partition, Lemma~\ref{lem:moment} gives
$h\log\E[e^{rJ_I/h}\mid\F_u]\le |I|\frac{r^2\kappa^2}{1-4r^2\kappa^2}$.

We now condition successively at the left endpoints, starting with the
last interval. Since all integrals over earlier intervals are
$\F_u$-measurable, each application of \eqref{eq:upper-cell} removes one
interval from the expectation. This gives
\begin{align*}
h\log\E[e^{rN_h^v/h}]
&\le \frac1r\sum_{I=[u,u+\delta]}\delta\sup_{x\in\R}\Psi(r^2a(u,x))+o(1)\\
&=\frac1r\int_0^1\sup_{x\in\R}\Psi(r^2a(t,x))\,dt+o(1).
\end{align*}
There are at most $1/\delta$ intervals of length $\delta$, so the
uniform $o(\delta)$ errors add up to $o(1)$. The two end intervals
contribute at most $O(\delta)$. Since $a$ is Lipschitz and $\Psi$ is uniformly continuous on
$[0,r^2\kappa]$, the function $t\mapsto\sup_{x\in\R}\Psi(r^2a(t,x))$ is continuous. The Riemann sums
therefore converge uniformly in $v$.
Combining this with \eqref{eq:upper-error-moment}, we conclude that
\[
\limsup_{h\downarrow0}h\log\E[e^{N_h(a)/h}]
\le \frac1{r^2}\int_0^1\sup_{x\in\R}\Psi(r^2a(t,x))\,dt.
\]
Finally, let $r\downarrow1$. Since $\Psi$ is uniformly continuous on
compact subintervals of $[0,1/2)$, we obtain the desired upper bound.
This completes the proof of Theorem~\ref{thm:nonlinear}.

\section*{Acknowledgments}
This research was supported by the Israel Science Foundation (ISF), grant 305/25.

\end{document}